\documentclass[12pt,a4paper]{article}

\newcommand{\F}{{F_{_0}(P)}}

\title{The Dimensions of the Symmetry Types of Polyhedra with Reflection Groups}

\author{\begin{tabular}{ccccc}
          {\small M. Rostami} & $\,$ & {\small Henrique F. da Cruz} & $\,$ & {\small Ilda I. Rodrigues} \\
          {\scriptsize rostami@ubi.pt} & $\,$ & {\scriptsize hcruz@ubi.pt} & $\,$ & {\scriptsize ilda@ubi.pt}
        \end{tabular}\\ {\small Universidade da Beira Interior}\\ {\small Covilh\~a - Portugal} }

\date{\today}

\newtheorem{lema}{Lemma}[section]
\newtheorem{teo}{Theorem}[section]

\newtheorem{obs}{Remark}[section]
\newtheorem{defi}{Definition}[section]

\usepackage{mathrsfs}
\usepackage{indentfirst}
\usepackage{setspace}

\usepackage[latin1]{inputenc}
\usepackage{fontenc} % fontes para os caracteres em portugues
\usepackage{amssymb}
\usepackage[dvips]{graphicx}
\usepackage{amsmath,a4}
\usepackage[dvips]{graphicx}
\usepackage{amsfonts}
\usepackage{latexsym}
\usepackage{enumerate}
\usepackage{fancyhdr}
\usepackage{makeidx}

\usepackage{epsfig}

\usepackage{stackrel}

\usepackage{multirow}
\usepackage{pstricks}
\usepackage{pst-node}
\usepackage{tikz}
\usepackage{pgfplots}

\begin{document}

\maketitle

\begin{abstract}
{\small{Let $P$ and $Q$ be convex polyhedra in $\mathbb{E}^3$ with  face lattices $F(P)$ and $F(Q)$ and symmetry groups $G(P)$ and $G(Q)$, respectively.  Then, $P$ and $Q$ are called face equivalent if there is a lattice isomorphism between  $F(P)$ and $F(Q)$; $P$ and $Q$ are called symmetry equivalent if the action of $G(P)$ on $F(P)$ is equivalent to the action of $G(Q)$ on $F(Q)$.  It is well known that the set $[P]$ of all polyhedra which are face equivalent to $P$ has the structure of a manifold of dimension $e-1$, up to similarities, where $e=e(P)$ is the number of edges of $P$. This is a consequence of the Steinitz's classical Theorem.  We give a new proof of this fact. The symmetry type  of $P$ denoted by $\langle P \rangle$ is the set of all polyhedron $Q$ symmetry equivalent to $P$. We show that for polyhedra with symmetry group $G(P)$ a reflection group the dimension of this manifold is $\epsilon-1$ where $\epsilon$ is the number of edge orbits of $P$ under and the action of $G(P)$ on $F(P)$.}}
\end{abstract}

\vspace*{0,5cm}

\vspace*{0,25cm}

\section{Introduction}

\indent Throughout this paper, the term ``polyhedron''  refers only to convex polyhedron in Euclidean 3-space $\mathbb{E}^3$ (convex 3-polytope). For any polyhedron $P$, we denote the set of all vertices, edges and faces of $P$ by $F_0(P)$, $F_1(P)$ and $F_2(P)$ respectively. These sets together with the empty set $\emptyset$ and $P$ itself, form a lattice $F(P)$ under inclusion, with $P$ as maximum and $\emptyset$ as minimum. The lattice $F(P)$ of $P$ is determined by the {\it incidence matrix} $M_P$ of $P$, where we label the vertices by $v_1,\ldots, v_r$ and the polygonal faces by $A_1,\ldots,A_s$, say, and define $M_P$ to be the $r \times s$ matrix whose $(i,j)th$ element $m_{i,j}$ is 1 if $v_i \in A_j$ and 0 otherwise. The number $\mu=\mu(M_P)=\mu(P)$, of nonzero elements of $M_P$ is called the {\it multiplicity} of $P$. Since each edge determines four incidences of the adjacent faces and vertices and each incidence corresponds to two edges, we have $\mu=\mu(P)=2e$, where $e=e(P)$ is the number of edges of $P$.

\indent Two polyhedra $P$ and $Q$ are face equivalent (also called combinatorially equivalent), written $P \approx Q$, if there is a lattice isomorphism from $F(P)$ to $F(Q)$. The set $[P]$ of all $Q \approx P$ is the \emph{face type} or, as often called, the \emph{isomorphism type} of $P$. Since $P$ is the convex hull of the set $F_0(P)$ of its vertices and since $P$ determines $\F$, we can identify $P$ with the point  $(v_1,\ldots,v_r)$ in $\mathbb{E}^{3r}=(\mathbb{E}^3)^r$, where $\F = \{ v_1,\ldots, v_r \} $. Thus we may identify a neighborhood of $P$ in $[P]$ with a neighborhood of $(v_1,\ldots, v_r)$ in $\mathbb{E}^{3r}$, since any $P'\approx P$ sufficiently close to $P$ may have its vertices labeled $v'_1,\ldots v'_r$ in a unique way so that $v'_i$ is close to $v_i$, $i=1,\ldots,r$.

A graph  $\mathscr{G}=(V,E)$   with vertex set  $V$  and edge set $E$  is said to be $3$-connected if any two vertices are connected by three internally disjoint paths (each pair of paths have only the two vertices in common); $\mathscr{G}$  is planar if it can be embedded in the plane such that no two edges intersect internally.

Now, for  any polyhedron $P$ each  edge has exactly two vertices. Therefore we can define the  $1$-skeleton, or edge graph, $\mathscr{G}_P$ of any  polyhedron  $P$ as  follows. The vertices of $\mathscr{G}_P$, are vertices of   $P$, and two vertices are connected by an edge in graph $\mathscr{G}_P$ if they are vertices of the same edge of $P$. This planar graph which is called ``\emph{Schlegel diagram}'' of $P$, can be obtained by projecting $P$ on one of its faces (outside face).

Due to Steinitz's famous Theorem called ``Fundamental theorem of convex types'' a graph $\mathscr{G}_P=(V,E)$  is the edge graph of  a polyhedron $P$ if and only if it is planar and $3$-connected. Beyond this remarkable result, it is interesting that, Steinitz discovered some important facts about the topological structure of $[P]$ which is called the realization  space  of $P$ (see Section 2 for the relevant definitions).

Indeed, a careful inspection of the proof of Steinitz's Theorem shows that $[P]$  is smooth manifold with dimension $\dim[P]=e-1$, up to similarities, where $e=e(P)$ is the number of edges of $P$ .

Perhaps the naturality of the above result of Steinitz was the main reason for other mathematicians to believe that the realization  space of polytopes, in dimensions higher than three, might also be a smooth manifolds. Theorem, on page 18 of \cite{Rob}, asserts that $[P]$ is manifold in general case, where $P$ is an $n$-dimensional convex polytope.  But a striking result of Mnev \cite{RG,Mnev} shows that, from the topological point of view, the realization  space of a convex polytope may be arbitrarily complicated. The argument given in \cite{Rob} does not take full account of coplanarity conditions on the intersections of the faces, but works for $n=3$ (see Section 3).

The structure of this note is as follows. In Section 2 it is introduced the explicit definition of the {\it realization space} of a polyhedron in terms of a relevant metric. Then, we give a proof that the realization space is in fact a smooth manifold with a specific dimension. The proof, which is independent of implicit proof of Steinitz himself, is illustrated with  an example adopted from (\cite{Rob}, p. 21).

%In section 3 we consider the polyhedra with symmetries generated by reflections.

In Section 3, for any polyhedron $P$ we define $\langle P \rangle$ the {\it symmetry type} of $P$, together with giving some basic facts about the actions of transformation groups on manifolds. Then, we recall the {\it stratification} of $\mathbb{E}^3$  and the decomposition of the manifold  of the space of all polyhedra into {\it strata} of symmetry types under the action of the point groups.

In Section 4, we first consider finite point groups generated by reflections. Then, with the help of the familiar notation of {\it fundamental region} of a reflection group, we prove our main theorem  which relates the dimension of the symmetry type of polyhedron $P$ to the number of its edge orbits, under the action of $G(P)$, the symmetry group of $P$ on its face lattice $F(P)$. This is the proof of Deicke's conjecture in  \cite{Rob} as well, for polyhedra having reflection symmetry groups. For polyhedra with rotation groups we refer to \cite{Ros1}, and for the basic properties of polyhedra with symmetry we refer to Robertson's book \cite{Rob}.

%This paper contains some parts of the first author unpublished doctoral thesis written under the supervision of Professor S. A. Robertson (University %od Southampton, U.K., 1987).

%We give a proof adopted from \cite{Rob} which is independent of "implicit proof" of Steinitz himself. Then we consider the symmetry part of this fact. For any polyhedron $P$ let
% $\langle P \rangle$  be the symmetry type of  $P$. We look to topological structure of $\langle P \rangle$ as well and prove that if $G(P)$ is a reflection group and the action of $G(P)$ on $F(P)$ has  $\epsilon$ edge orbits then the dimension of symmetry type of $P$ , $\dim \langle P \rangle= \epsilon -1$. This is the proof of Deicke's conjecture \cite{Rob} for reflection groups. For polyhedra with rotation groups we refer the reader to \cite{Ros1}.

\section{The realization space of a polyhedron}

\indent Our purpose in this Section is to define realization space and determine its topological structure. But first we present some basic definitions.

Let $\left( \mathbb{E}^n , \, d \right)$ be the $n$-dimensional Euclidean space with induced norm
%\vspace*{-.2cm}
$$
d(x,y) \, = \, \|x-y \|, \, x,y\in \mathbb{E}^n \, .
$$

%\vspace*{-.2cm}
Denote by $\mathbb{K}^n$ hyperspace of all nonempty compact convex subsets (convex bodies) of $\mathbb{E}^n$. Then, $\mathbb{K}^n$ is endowed with the following familiar metric, called \emph{Hausdorff metric}.

\begin{defi}
For $P\in \mathbb{K}^n$ and $\varepsilon >0$, let
%\vspace*{-.2cm}
$$
\mathcal{U}_{\varepsilon}(P) \, = \, \left\{ x\in \mathbb{E}^n \, \mid \, d(x,P) < \varepsilon \right\} \, ,
$$

%\vspace*{-.2cm}
\noindent where $d(x,P) \, = \, \stackrel[p\in P]{}{\inf} \left\{ d(x,p)=\|x-p \| \right\}$ .\\

%\vspace*{0,2cm}
Now, for $P, \, P' \in \mathbb{K}^n$, let $\rho(P,P') \, = \, \inf  \left\{ \varepsilon \, \mid \, P' \subseteq \, \, \mathcal{U}_{\varepsilon}(P) \right\}$. Then,
%\vspace*{-.2cm}
$$\,
\,
\begin{array}{rcl}
  d_{\mathcal{H}}(P,P') & = & \max  \left\{ \rho(P,P') \, , \, \rho(P',P) \right\} \\[4pt]
                        & = & \inf \left\{ \varepsilon >0 \, \mid \, P \subseteq \, \, \mathcal{U}_{\varepsilon}(P') \text{ and } P' \subseteq \, \, \mathcal{U}_{\varepsilon}(P) \right\} \, ,
\end{array}
$$

%\vspace*{-.2cm}
\noindent is a metric on $\mathbb{K}^n$ which is called Hausdorff metric (distance).
\end{defi}

We denote by $\mathscr{P}$ the set of all polyhedra in $\mathbb{E}^3$ as a topological subspace of the metric space $\left( \mathbb{K}^3, \, d_{\mathcal{H}} \right)$.

\begin{defi}
Let $P$ be a polyhedron (convex 3-polytope). Denote by $[P]$ the face type of $P$, the set of all polyhedra $Q$ face equivalent to $P$.
%\vspace*{-.2cm}
$$
[P] \, = \, \left\{ Q \, \mid \, Q {\text{ is polyhedron and }} Q\approx P \right\} \, .
$$

%\vspace*{-.2cm}
\noindent Then, $[P]$ together with its natural subspace topology induced by Hausdorff metric $d_{\mathcal{H}}$ is called the realization space of $P$.
\end{defi}

Alternatively, any polyhedron $P$ determines, and it is determined by, the set $F_0(P) = \{ v_1,\ldots,v_r \}$ of its vertices. Since we can identify $P$ with the point $(v_1,\ldots,v_r)\in \left( \mathbb{E}^3 \right)^r$, $[P]$ can be interpreted as a topological subspace of $\mathbb{E}^{3r}$ with its topology induced by the vertices.

Given $\varepsilon >0$, sufficiently small, there is an open neighborhood $V_{\varepsilon}(P)$ of $P$ in $[P]$ such that, for all $P' \in V_{\varepsilon}(P)$ we have $F_0(P') = \{ v_1^{'},\ldots,v_r^{'} \}$, where $\| v_i - v_i^{'} \| < \varepsilon$, $i=1, \ldots , r$.

Therefore, $[P]$ can be topologized locally by $V_{\varepsilon}(P)$ neighborhood of $P$. Thus $V_{\varepsilon}(P)$, in fact, is the open space of the small perturbations of the vertices of $P$.

\vspace*{0,2cm}
Now, let $\Pi_j$ be the plane containing $A_j$. Then, we may suppose without loss of generality that for all $j=1,\ldots,s$, the origin $O$ does not lie in $\Pi_j$. Otherwise translate  $P$ to ensure this condition. Then, for all $j$ there is a unique $a_j \in \mathbb{E}^3$ such that $\Pi_j$ is given by the equation $\langle x,a_j\rangle=1, x\in \mathbb{E}^3$. Thus $\Pi_j$ is given by $a_j$ and $P$ itself by the point
%\vspace*{-.1cm}
$$(v_1,\ldots,v_r,a_1,\ldots,a_s)\in \mathbb{E}^{3r}\times \mathbb{E}^{3s}=\mathbb{E}^{3(r+s)} \, ,$$

%\vspace*{-.1cm}
\noindent where $\langle v_i,a_j\rangle=1$ for all $i,j$ such that $m_{i,j}=1$. Now let $P'$ be the polyhedron with $P\approx P'$, having vertices $v'_1,\ldots,v'_r$ and faces given by $a'_1,\ldots,a'_s$ where the plane $\Pi'_j$ of $A'_j$ has equation $\langle x,a'_j\rangle= 1$ and again $\langle v'_i,a'_j \rangle=1$ if $m_{ij}^{'}=1$ (with obvious notation).
Suppose that $P'$ is close to $P$, so that
%\vspace*{-.3cm}
$$v'_i=v_i+\xi_i,\,\,\,\,a'_j=a_j+\eta_j \, ,$$

%\vspace*{-.3cm}
\noindent for some $\xi_j,\eta_j \in \mathbb{E}^3$ with $||\xi_i||$ and $||\eta_j||$ small ($i=1,\ldots,r;j=1,\ldots,s$). Then,
%\vspace*{-.2cm}
\begin{equation}\label{eq1}
\bf{\Phi_{[i,j]}(\xi,\eta):=\langle v_i,\eta_j\rangle+\langle \xi_i,a_j\rangle+\langle\xi_i,\eta_j\rangle=0},
\end{equation}

%\vspace*{-.3cm}
\noindent if $m_{ij}=1$ (and hence $m_{ij}^{'}=1$).

\vspace*{0,2cm}
Conversely, for any sufficiently small $\xi \in \mathbb{E}^{3r}$ and $\eta \in \mathbb{E}^{3s}$, the point
%\vspace*{-.3cm}
$$(\nu_1+\xi_1,\ldots,\nu_r+\xi_r,a_1+\eta_1,\ldots,a_s+\eta_s)$$

%\vspace*{-.3cm}
\noindent represents a unique polyhedron $P'\approx P$, provided equations (\ref{eq1}) hold. Let us now consider the polynomial map $\Phi: \mathbb{E}^{3r}\times \mathbb{E}^{3s}\rightarrow \mathbb{E}^\mu$ given by  $\Phi(\xi,\eta)=(y_1,\ldots,y_\mu)$ where the pairs $(i,j)$ with $m_{ij}=1$ are arranged in lexicographical order, and if $(i,j)$ is the $[i,j]th$ such pair then, $y_{[i,j]}=\Phi_{[i,j]}(\xi,\eta)$ is defined as in (\ref{eq1}). We may for convenience suppress component suffices of $\xi_i$ and $\eta_j$, and write
%\vspace*{-.3cm}
$$\frac{\partial y_{[i,j]}}{\partial \xi_i}=a_j+\eta_j,\,\,\,\,\frac{\partial y_{[i,j]}}{\partial \eta_j}=\nu_i+\xi_i \, .$$

%\vspace*{-.3cm}
Thus the Jacobian matrix $J_{\Phi}(0,0)$, has order  $\mu\times(3(r+s))$ with $a_j$ in the $[i,j]$th row and the $i$th triplet of columns, and $\nu_i$ in the $[i,j]$th row and the $(r+j)$th triplet of columns. All the other elements of $J_{\Phi}(0,0)$ are $0$.

A simple example may helps to clarify these remarks. Let $P$ be the polyhedron shown, with its accompanying Schlegel diagram in the following figure.
%UMA FIGURA
{ \begin{figure}[!htb] \centering
\hspace*{.75cm} \includegraphics[scale=.2]{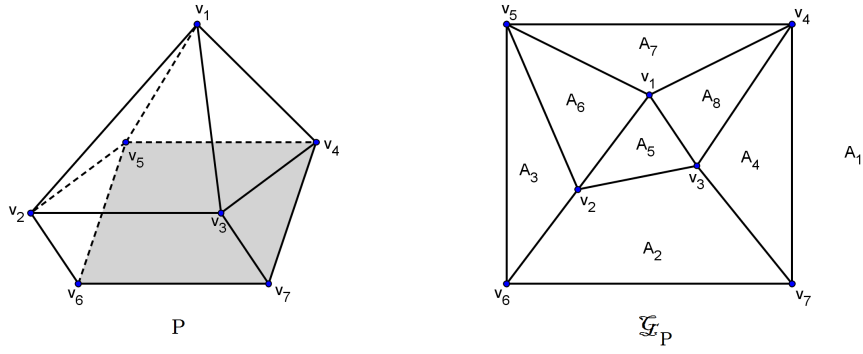}
\vspace*{-0,4cm}\caption{{\small{Schlegel diagram of $P$ with ``outside face'' $A_1$}}}
{\label{im1}}
\end{figure}}

\newpage
In this example, $r=7,s=8$ and $\mu=26$. Notice that $3(r+s)=45>26=\mu$. In fact, for any polyhedron, $3(r+s)=3(e+2)=3e+6=\mu+e+2>\mu$, by Euler's Theorem. Thus $\mu$ is the largest value of the rank of $J_{\Phi}(0,0)$. In the example, $J_{\Phi}(0,0)$ has 26 rows and 45 columns, shown below in the truncated $26\times 15$ form.

\begin{center} {\scriptsize{\begin{tabular}{|c|c|c|c|c|c|c|c|c|c|c|c|c|c|c|c|}
  \hline
    & 1 & 2 & 3 & 4 & 5 & 6 & 7  & 1 & 2 & 3 & 4 & 5 & 6 & 7   & 8  \\
    \hline
 15 & $a_5$ &   &   &   &   &   &     &   &   &   & $\nu_1$   &  &  &  & \\
  \hline
  16 & $a_6$ &   &   &   &   &   &     &   &   &   &   & $\nu_1$ &   & & \\
  \hline
  17 & $a_7$ &   &   &   &   &   &     &   &   &   &   &   &  & $\nu_1$ & \\
  \hline
  18 & $a_8$ &   &   &   &   &   &     &   &   &   &   &   &   &   &$\nu_1$ \\
  \hline
  22 &   & $a_2$ &   &   &   &   &     &  & $\nu_2$  &   &   &   &   &   & \\
  \hline
  23 &  & $a_3$ &   &   &   &   &      &         &  &  $\nu_2$ &   &   &   &   & \\
  \hline
  25 &   & $a_5$ &   &   &   &   &     &   &   &  &  & $\nu_2$  &   &   & \\
  \hline
  26 &   & $a_6$ &   &   &   &   &     &   &   &   &   &  & $\nu_2$  &   & \\
  \hline
  32 &   &   & $a_2$ &   &   &   &     &  & $\nu_3$  &   &   &   &   &   & \\
  \hline
  34 &   &   & $a_4$ &   &   &   &     &   &  &   & $\nu_3$ &   &   &   & \\
  \hline
  35 &   &   & $a_5$ &   &   &   &     &   &   &   &  & $\nu_3$  &   &   & \\
  \hline
  38 &   &   & $a_8$ &   &   &   &    &   &   &   &   &   &   &   & $\nu_3$\\
  \hline
  41 &   &   &   & $a_1$ &   &   &   & $\nu_4$   &   &   &   &   &   &  & \\
  \hline
  44 &   &   &   & $a_4$ &   &   &   &     &   &  & $\nu_4$  &   &   &   & \\
  \hline
  47 &   &   &   & $a_7$ &   &   &   &     &   &   &   &   &  & $\nu_4$ & \\
  \hline
  48 &   &   &   & $a_8$ &   &   &   &     &   &   &   &   &   &   &$\nu_4$\\
  \hline
  51 &   &   &   &   & $a_1$ &   &   & $\nu_5$   &   &   &   &   &   &   & \\
  \hline
  53 &   &   &   &  &  $a_3$ &   &   &     &  & $\nu_5$  &   &   &   &   & \\
  \hline
  56 &   &   &   &   & $a_6$ &   &   &      &   &   &   &  &  $\nu_5$ &   & \\
  \hline
  57 &   &   &   &   & $a_7$ &   &   &     &   &   &   &   &  & $\nu_5$  & \\
  \hline
  61 &   &   &   &   &   & $a_1$ &   &  $\nu_6$  &   &   &   &   &   &   &\\
  \hline
  62 &   &   &   &   &   & $a_2$ &     &  &  $\nu_6$ &   &   &   &   &   &\\
  \hline
  63 &   &   &   &   &   &$ a_3$ &     &   &  & $\nu_6$  &   &   &   &   &\\
  \hline
  71 &   &   &   &   &   &   & $a_1$ &  $\nu_7$ &   &   &   &   &   &   & \\
  \hline
  72 &   &   &   &   &   &   & $a_2$ &   & $\nu_7 $  &   &   &   &   &   & \\
  \hline
  74 &   &   &   &   &   &   & $a_4$ &    &   &  &   $\nu_7$ &   &   &   &\\
  \hline
\end{tabular}}} \end{center}

\begin{center}
%Table 1:
[truncated $J_{\Phi}$]
\end{center}

%\vspace*{.25cm}

Recall that a similarity is a map
\vspace*{-.3cm}
$$ s:\mathbb{E}^3 \rightarrow \mathbb{E}^3$$

\vspace*{-.3cm}
\noindent such that for some real $r>0$, and all $x,y\in \mathbb{E}^3$
\vspace*{-.3cm}
$$
rd(x,y)=d(s(x), s(y)) \, .
$$
%  where $d$ is usual metric in $\mathbb{E}^3$.

\vspace*{-.3cm}
The set of all similarities is a group which is denoted by $Sim(3)$, and called similarity group of $\mathbb{E}^3$.

The map $\varphi: Sim(3)\rightarrow \mathbb{R}_\ast$ defined by $\varphi(f)=r$, where $r$ is as above and $\mathbb{R}_\ast$ is multiplicative group of positive real numbers, is a group homomorphism.
It is well known that the kernel of this homomorphism is the isometry group or the Euclidian group $\mathbb{E}(3)$. We may identify $Sim(3)$ as $\mathbb{E}(3) \times \mathbb{R}_\ast$. Since the manifold $\mathbb{E}(3)$ has dimension $6$, the dimension of $Sim(3)$ is equal to $7$.

\vspace*{.2cm}
Now we are ready to prove our main theorem in this Section:
\begin{teo} Let $P$ be a convex 3-polyhedron. Then, $[P]/Sim(3)$, the realization space of $P$ modulo similarities, is a manifold of dimension $\dim \left( [P]/Sim(3) \right) \, = \, e-1$, where $e=e(P)$ the number of edges of $P$.
\end{teo}
\vspace*{-.1cm} {\bf Proof:} To show, equivalently, that for any $P$, the face type $[P]$ is a manifold of dimension $e+6$, it is enough to show that $J_{\Phi}(0,0)$ has rank $\mu$. The result follows by the implicit function Theorem, since $e-1=3(r+s)-\mu-7$.

We want to show that the rows of $J_{\Phi}(0,0)$ are linearly independent. Suppose then, that for some real number $\alpha_{[i,j]}$, we have
\vspace*{-.1cm}
\begin{equation}\label{eq2}
\sum_j \alpha_{[i,j]}a_j=\sum_i \alpha_{[i,j]}\nu_i=0 \, ,
\end{equation}

\vspace*{-.3cm}
\noindent for each $i=1,\ldots,r$ and each $j=1,\ldots,s$. We know that for any $i$ and any three values of $j_1,j_2,j_3$ of $j$ with $m_{ij}=1$, $a_{j_1},a_{j_2},a_{j_3}$ are linearly independent. Hence for each $i$, the equations $\sum  \alpha_{[i,j]}a_j=0$ have solution space of dimension $s_i-3$ where $s_i$  is the number of the values of $j$ with $m_{i,j}=1$, and we can express any three of the  numbers  $\alpha_{[i,j]}$ as linear functions  of the remaining $s_i-3$. Hence the $\mu$ numbers $\alpha_{[i,j]}$ are expressible as linear functions of $\sum_{i=1}^r(s_i-3)=\mu-3r$ independent variables. In particular, for each $i$ for which $s_i=3$, all the numbers $\alpha_{[i,j]}$ are 0.

Now consider the equations $\sum \alpha_{[i,j]}\nu_i=0$. Again, for each $i$, any three of the vertices $v_i$ are linearly independent. So again we can express any three of the numbers $\alpha_{[i,j]}$ as linear combinations of the remaining $r_j-3$. So the numbers $\alpha_{[i,j]}$ are expressible as linear functions of $\sum_{i=1}^s(v_j-3)=\mu-3s$ independent variables.

But the labels $\alpha_{[i,j]}$ for each $j$ are already expressed as linear functions of $\mu-3r$ independent variables, as described above. Hence the numbers $\alpha_{[i,j]}$ are over determined by the two sets of equations in (\ref{eq2}), since the solution space of the equations in (\ref{eq2}) has dimension $\mu-3(r+s)$. But  $\mu-3(r+s)<0$ by Euler formula. It follows that $\alpha_{[i,j]}=0$ for all $i,j$ where $m_{i,j}=1$.\hfill$\square$

%$\mu-3r-3s<0$, it follows that $\alpha_{[i,j]}=0$ for all $i,j$ where $m_{i,j}=1$.$\diamondsuit$

\vspace*{.5cm}
Referring to the above example, we may express each of, say  $\alpha_{[1,6]}$, $\alpha_{[1,7]}$ and $\alpha_{[1,8]}$ as linear function of $\alpha_{[15]}$ that is as a multiple of $\alpha_{[15]}$. Thus $\alpha_{[16]}=\lambda_{16}\alpha_{[15]}$, $\alpha_{[17]}=\lambda_{17}\alpha_{[15]}$, $\alpha_{[18]}=\lambda_{18}\alpha_{[15]}$. Likewise, taking the first number $\alpha_{[i,j]}$ in each block parameter where possible we have  $\alpha_{[23]}=\lambda_{23}\alpha_{[22]}$, $\alpha_{[25]}=\lambda_{25}\alpha_{[22]}$ and $\alpha_{[26]}=\lambda_{26}\alpha_{[22]}$, and so for each of the values $i=1,2,3,4,5$. For $i=6$ and for $i=7$, $s_i=3$ and immediately we find that $\alpha_{[6i]}=\alpha_{[7j]}=0$ for each $j$ with $m_{6j}=1$ or $m_{7j}=1$.

Now repeat this procedure with numbers $\alpha_{[i,j]}$ with $j$ fixed. Thus  $\alpha_{[51]}=\mu_{51}\alpha_{[41]}$, $\alpha_{[61]}=\mu_{61}\alpha_{[41]}$, and $\alpha_{[71]}=\mu_{71}\alpha_{[41]}$, and likewise, for $j=2$. For $j=3,4,5,6,7$ and $8$, only three values of $i$  corresponds to each values of $j$, because the associated faces $A_j$ are triangles. So all $\alpha_{[ij]}=0$ for all $i,j$ with $m_{i,j}=1$, as in general case.

\begin{obs}
{\label{obs1}}
%\textbf{Remark 2.1}
{\rm It is worth mentioning that Richter-Gebert in (\cite{RG}, Section 13.3) by fixing a suitable affine basis, proves that the realization space of a polyhedron $P$ (denoted there by $\mathcal{R}(P)$), that is the space of coordinatization for the combinatorial type of $P$ with $e(P)=e$ edges is a smooth open manifold of dimension $e-6$, module the natural action of $12$-dimensional affine transformation group.

The realization space $\mathcal{R}(P)$ is understood as a subspace of $\mathbb{E}^{3n}$ by identifying the $3n$ coordinates of the $n$ vertices of $P$ with points in $\mathbb{E}^{3n}$. It is described by the set of all solutions of a collection of polynomial equations and inequalities with integer coefficients. Such sets are called a simple semi-algebraic variety.

Furthermore, by fixing an affine basis in the definition of the realization space, one makes sure that the ``reflection'' (mirror images in Steinitz's proof) do not create a second component of the realization space. Therefore, $\mathcal{R}(P)$ is indeed (path connected) \emph{contractible} and has \emph{Steinitz's isotopy property} (\cite{Stei}, Section 69), i. e., any two realizations $P_1$ and $P_2$ of $P$ can be continuously deformed into each other while maintaining the same structure throughout.

However, our approach to study the face type manifold $[P]$ of a polyhedron $P$, is different. Although the realization spaces such as $[P]$ are usually (for example as above) defined modulo affine or Euclidian groups, in this work we consider $[P]/Sim(3)$, $[P]$ modulo similarity group. The topology of the realization space $[P]$ of polyhedron $P$ is induced by Hausdorff metric, and since the action $Sim(3)$ on $[P]$ is not (fixed-point) free (see Section 3 for definition), the manifold $[P]/Sim(3)$ is not contractible, and clearly does not satisfy Steinitz's isotopy property.}
\end{obs}

\begin{obs}
{\rm Let $P$ be a polyhedron with $f_0(P)=r$ vertices, $f_1(P)=e$ edges and $f_2(P)=s$ faces. Consider the face type $[P]$ of $P$. In (\cite{Rob}, p. 75) the dimension of this manifold is given by
%\vspace*{-.2cm}
$$\dim [P]=3(r+s)-\mu(P) \, .$$

%where $e$ is the number of vertices, $f$ the number of faces of $P$, and $\mu (P)=2e$ is the multiplicity of $P$ (see \cite{Stei}, p. 75).
%\vspace*{-.3cm}
An intuitive derivation of the above formula may be given as follows. If the vertices and faces of $P$ were allowed to move independently in $\mathbb{E}^3$ then, they would have $3(r+s)$ \emph{degrees of freedom} (see page 11 for definition). But they are not independent. In fact for each incidence relation $(v,F)\in F_0(P)\times F_2(P)$ with $v \in F$, the whole space loses one degree of freedom. Hence we have
%But any $k$ vertices to specie a $k$-gonal face $k\geq 3$, they loose $k-3$ degrees of freedom. Thus the whole configuration has
%\vspace*{-.3cm}
$$\dim[P] \, = \, 3(r+s)-\mu(P)\,.$$

%\vspace*{-.3cm}
%\noindent degrees of freedom; the sum being taken over all faces of $P$. Now if we denote by  $x_k$ the number of $k$-gonal faces and by $f$ the total number of all faces of P, then since each edge is contained in exactly two faces of $P$, we have $\sum x_k =2e$ and
%\vspace*{-.3cm}
%$$\dim [P] \, = \, \sum_{k\geq 3} (k-3) \, = \, 3r-\sum kx_k + 3s \, = \, 3(r+s)-\mu(P) \, .$$
%$$\begin{array}{rcl}
%  \dim [P] & = & \sum_{k\geq 3} (k-3) \\[10pt]
%    & = & 3r-\sum kx_k + 3s \\[10pt]
%    & = & 3(r+s)-\mu(P) \, .
%\end{array}$$

%\vspace*{-.3cm}
Note that, since the similarity Lie group $Sim(3)$ has dimension 7, by Euler's formula
\vspace*{-.1cm}
$$n-e+f=2 \: , $$

\vspace*{-.1cm}
\noindent we have
\vspace*{-.1cm}
$$
\begin{array}{rcl}
  \dim \left( [P]/Sim(3) \right) & = & 3(r+s)-\mu(P)-7 \\
  \, & = &  3(r+s)-2e-7 \, = \, e-1\, (\text{\small{\emph{module similarities}}}) \, .
\end{array}
$$}
\end{obs}

\section{Symmetry type of a polyhedron and stratifications}

{\bf Transformation groups}

\vspace*{0,2cm}
\indent We start this Section by introducing some basic definitions and standard results from the theory of Lie groups acting on smooth manifolds and refer the reader to \cite{Dei} and \cite{Kaw} for more details.

\vspace*{0,2cm}
Let $G$ be a Lie group, with identity element $e$, and $M$ a smooth manifold. A smooth action of $G$ on $M$ is a $C^{\infty}$ mapping:
\vspace*{-0,2cm}
$$\phi: G\times M \rightarrow M \: , \:\: \phi(g,x)=\phi_g(x)=g(x) \, ,$$
such that
$$
\begin{array}{rcll}
  e(x) & = & x, & x\in M \, ,\\
  (g_1g_2)(x) & = & g_1(g_2(x)), & g_1, g_2 \in G \, , \: x\in M \, .
\end{array}
$$

In this case we say that $M$ is a $G$-manifold. For any $x\in M$ the subgroup
\vspace*{-.2cm}
$$G_x=\{g \in G \mid g(x)=x\} \, ,$$
of $G$ is called the \emph{stabilizer} or \emph{isotropy} subgroup at $x$.

The set $ \, G(x)=\{ g(x) \mid g\in G \} \,$ is called the \emph{$G$-orbit} of $x$. The \emph{orbit space} of the action of $G$ on $M$ is the space $M/G$, the space of all $G$-orbits endowed with the quotient topology given by canonical projection
\vspace*{-.2cm}
$$
\begin{array}{ccccl}
  \pi & : & M & \rightarrow & M/G\\
   &  & x & \mapsto & G(x)
\end{array} \, ,
$$
\noindent and the differentiable structure of $M/G$ is induced by the same structure of $M$.
% estava   differentiable              e subtituimos por SAME

The action is called \emph{free} if, for each $x\in M$, $ \: G_x=\{ e \} \, $ .

If a Lie group $G$ acts on a smooth manifold $M$ via $\phi$, we call $(M,G):=(M,G,\phi)$ a transformation group, and $M$ is said to be a $G$-manifold.

%\vspace*{0,2cm}
Now, for each polyhedron $P$, a symmetry of  $P$ is a rigid  transformation  (or self isometry) $f:\mathbb{E}^3\rightarrow \mathbb{E}^3$  such that $f(P)=P$. Any such symmetry maps vertices to vertices, edges to edges and faces to faces and preserves inclusions (incidences). Hence any symmetry induces an  automorphism on $F(P)$. The set $G(P)$  of all symmetries of $P$ is a finite subgroup of the Euclidean group  $\mathbb{E}(3)$ acting on $F(P)$ as a group of automorphisms. We may assume that the centroid is $O$, so $G(P)$ is a finite subgroup of the orthogonal group $\mathcal{O}(3)$. If a finite subgroup $G$ of $\mathcal{O}(3)$ is the symmetry group of a convex polyhedron $P$, we also call $P$ a $G$-polyhedron.

%Then, we call $P$ a $G$-polyhedron.

\begin{defi}
Two polyhedra $P$ and $Q$ are \emph{symmetry equivalent}, and write $P\cong Q$, if there is an isomorphism
%\vspace*{-.3cm}
$$\lambda: F(P)\rightarrow F(Q)$$

%\vspace*{-.2cm}
\noindent of the face lattices and some isometry $f:\mathbb{E}^3 \rightarrow \mathbb{E}^3$ such that for all
$g \in G(P)$ and all $x \in F(P)$,
%\vspace*{-.1cm}
$$\lambda(gx)=(f g f^{-1})(\lambda(x)) \, .$$
\end{defi}

%\vspace*{-.3cm}
If we further assume that $P$ and $Q$ are both $G$-polyhedra, that is having the same (rather than conjugate) subgroups then, with the above condition, we say that $P$ and $Q$ are $G$-equivalent. Hence, in this case
%\vspace*{-.3cm}
$$\lambda(gx)=g\lambda(x) \, .$$

%\vspace*{-.5cm}
\begin{defi}
Let $P$ be a $G$-polyhedron. The symmetry type $\langle P \rangle$ of $P$ is defined by
%\vspace*{-.2cm}
$$\langle P \rangle = \{ Q \, | \, Q \text{ is $G$-polyhedron and $Q$ is $G$-equivalent to $P$}, \,\, Q\cong P \} \, . $$
\end{defi}

Now consider $\mathscr{P}$ the space of all convex polyhedra in $\mathbb{E}^3$. Since the subdivision of $\mathscr{P}$ into face types and symmetry types both respect the Euclidian similarities, it is convenient to look at the action of the similarity group $Sim(3)$ on the space $\mathscr{P}$ of all polyhedra in $\mathbb{E}^3$. This action partitions the quotient space
%\vspace*{-.2cm}
$$ \, \mathscr{S}:=\mathscr{P}/Sim(3) \, ,$$

\vspace*{-.2cm}
\noindent of similarity classes or ``shapes'' of polyhedra, into orbit types, where each orbit type consists of all those orbits on which the isotropy subgroups at any polyhedron in the orbit are conjugate.

Thus the symmetry types partitions each $[P]$ into mutually disjoint subsets refining the partitions of $\mathscr{S}$ into face types.

But the isotropy subgroup at $P$ is just the symmetry group $G(P)$ itself. It follows that the symmetry types are composed of components of the orbit types. The principal orbit type corresponds to the trivial isotropy subgroup, that is to say to the symmetry type of any polyhedron $Q$ in $[P]$ with trivial symmetry group $G(Q)=\{ e \}$. This type is open in $[P]/Sim(3)$ of dimension $e(P)-1$. All other symmetry types have lower dimensions.

As an example, let $P$ be a polyhedron combinatorially equivalent to cube. Then, $[P]/Sim(3)$ has dimension 11. The principal orbit type corresponds to the realization space of polyhedron $Q\approx P$ with trivial group, is an open and dense submanifold of $[P]/Sim(3)$ of dimension $e(Q)-1=11$.

Now, since each orbit type is a submanifold of $[P]$, by the Slice Theorem of transformation groups (see \cite{Kaw}, Th. 4.11, and \cite{Rob}, p.42), we have the following well known theorem.

%each symmetry type $\langle P \rangle$ of $P$ being, in effect, an orbit type of the action of the similarity group $Sim(3)$, the realization space $[P]$, is a smooth manifold.

\begin{teo} {\rm (in \cite{Rob}, p. 42 and \cite{Ros2})}  Let $P$ be a polyhedron in $\mathbb{E}^3$. Then, $\langle P \rangle$, the  symmetry type of $P$, is a smooth manifold.\end{teo}

In order to study the dimension of the symmetry type of $P$, $\dim \langle P \rangle$, we simply study those polyhedra $Q$ that lie in some neighborhood of $P$ in $\mathscr{S}$ and are symmetry equivalent to $P$. We can simplify our discussion by restricting our attention to those $Q$ whose symmetry group is not merely conjugate to $G(P)$ but $G(P)$ itself. Along this restriction we also factor out components that came from similarities. In this way, we may determine the value of  $\dim \langle P \rangle$. We first look at a simple example again.

Let $P$ be the right pyramid over a square (Figure \ref{image1}) with symmetry group $G(P)$ which is dihedral reflection group.

%Now consider $G(P)$ the symmetry group of $P$ which is dihedral reflection group.

Suppose we fix the group $G = G(P)$. Then, vertex $v_1$ can be chosen only on the axis of $G$. Therefore it has one degree of freedom (see page 11). Likewise $v_2$  must lie on reflection plane, hence has only two degrees of freedom. Having chosen  $v_2$, the vertices $v_3$,$v_4$ and $v_5$  which are on the same orbit of $v_2$ , have no degree of freedom at all, since they are determined by our choice of G and $v_2$ . Hence the vertices have a total of $1 + 2 = 3$ degrees of freedom. Similarly for the faces, each triangular face or equivalently the plane that contains it has only two degrees of freedom in the space of affine plane in $\mathbb{E}^3$, since each plane is invariant under a reflection element of $G$. But the square face has only one degree freedom, because it is orthogonal to the axis of rotation of $G$. Therefore the faces have just $2 + 1 = 3$ degrees of freedom. Of course the faces and vertices can not be chosen independently of one another. The incidence of $v_1$ with respect to any of the four triangular faces adjacent to it, determines the incidence of that vertex to the other three faces under the action of $G$. Hence $v_1$  has only one ``independent'' incidence. The vertices $v_2$, $v_3$, $v_4$ and $v_5$ are in the same $G$-orbit and each one is incident with two triangles and one square faces. Take one of them say $v_2$. There is a reflection which fixes $v_2$ and sends adjacent triangular faces each one to the other. Thus the number of independent multiplicity (to be defined later) of $P$ is $1 + 2 = 3$. Each such incidence relation in the form of the condition that a vertex lies in a particular face, reduces the dimension of the symmetry type by one. Now, if we take into account the fact that the center of $P$ can be chosen only on the fixed point set of $G$, which is one dimensional and considering also the dilation of $P$ which in each case reduces $\dim\langle P\rangle $ by one, we get
$$\dim\langle P\rangle \, = \, (1 + 2) + (1 + 2) - (1 + 2) - 2 \, = \, 1 \, .$$

\vspace*{.5cm}
%UMA FIGURA
{ \begin{figure}[!htb] \centering
\vspace*{-.5cm} \hspace*{.75cm} \includegraphics[scale=.25]{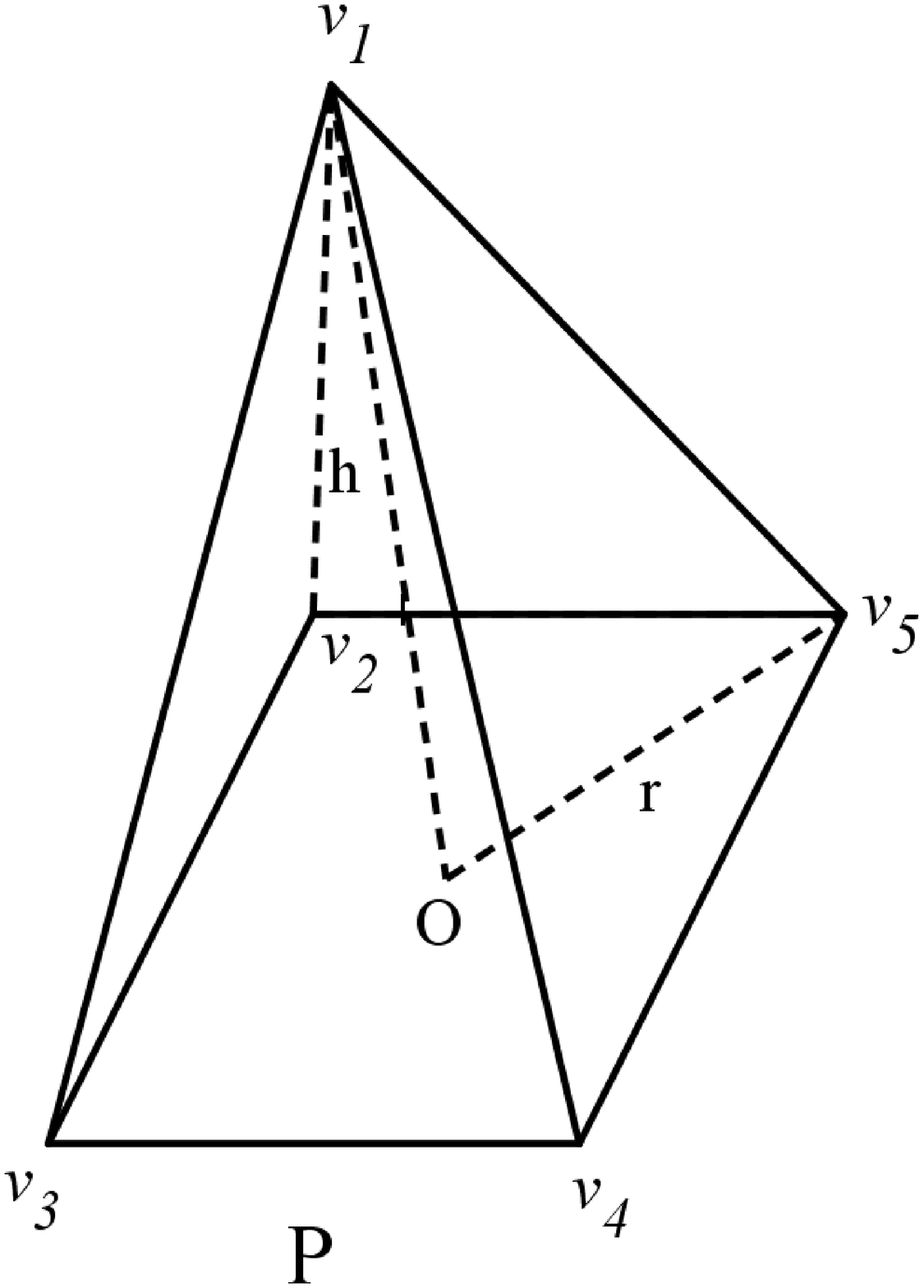}
\vspace*{-0,4cm}\caption{$\,$}
\label{image1}
\end{figure}}

\vspace*{.5cm}
In fact, in Figure \ref{image1} or in any right pyramid with a regular polygon as base, if we denote the height and radius of the base of P by $h$ and $r$, respectively and consider the ratio $\zeta=\frac{h}{r}$  the two such pyramids are similar if they have the same ratio $\zeta$. Therefore we can parameterize the symmetry type of $P$ by $\zeta$ with $0<\zeta$. Hence $\langle P \rangle$ has the structure of the open interval and $\dim \langle P\rangle  = 1$. Indeed the action of $G(P)$ on $F(P)$ has $\epsilon=2$ edge orbits. Hence
$\dim \langle P\rangle =\epsilon -1 = 1$.  Note that in the right pyramid with regular base the ratio $\zeta$ is similarity invariant and the symmetry type is a connected 1-manifold with boundary 0-dimensional symmetry types one for regular base and the other a segment (1-polytope). The following figure illustrates this idea.

\vspace*{.5cm}
\begin{center}
\begin{tikzpicture}[scale=0.5]
\filldraw [black]  (0,0) circle (2pt)
            (15,0) circle (2pt);
\draw (15,0) -- (0,0) node[above] {\small{$h\rightarrow 0$}};
\draw (15,0) -- (0,0) node[below] {\footnotesize{regular base}};
\draw (0,0) -- (15,0) node[above] {\small{$r\rightarrow 0$}};
\draw (0,0) -- (15,0) node[below] {\footnotesize{1-polytope}};
\end{tikzpicture}
\end{center}

\begin{center}
\vspace*{-1.1cm} {\footnotesize{right pyramid}}
\end{center}

%UMA FIGURA
\begin{center}
{ \begin{figure}[!htb] \centering
\vspace*{-.55cm} \includegraphics[scale=.075]{extra.eps}
\caption{$\,$}
\label{extra}
\end{figure}}
\end{center}

\vspace*{-1cm}
The idea of this example can be applied in general to find the dimension of the symmetry type of any polyhedron $P$.

%\vspace*{.75cm}
Let us denote by $F_0(P), F_1(P)$ and $F_2(P)$ the set of all vertices, edges and faces of $P$, respectively, with symmetry group $G=G(P)$.
%with orders $n$, $e$ and $f$ respectively.
%Let $G=G(P)$ be symmetry group of $P$.

\begin{defi}\label{TwoP}
Two ordered pairs $(v, F)$   and   $(v', F')$   in   $F_0(P)\times F_2(P)$   are  called  $G$-independent incidences or simply independent incidences, if and only if, there exists no $g \in G$   such   that  $g(v)=v'$    and $g(F)=F'$.  By  $\mu_*(P)$   we  mean  the  number  of  $G$-independent incidences $(v, F)$ where $v \in F_0(P)$,  $F \in F_2(P)$  and $v \in F$. Therefore $G(P)$ acts on the set of all such incident pairs $(v, F)$ with  $\mu_*(P)$  orbits of independent incidences or ``incident orbits''.
\end{defi}

For example let $P$ be a rhombic dodecahedron (Figure \ref{f2}) then, $\mu(P) = 2e = 48$ but  $\mu_*(P)=2$.
%UMA FIGURA
\vspace*{-.25cm}
{ \begin{figure}[!htb] \centering
\hspace*{.5cm} \includegraphics[scale=0.15]{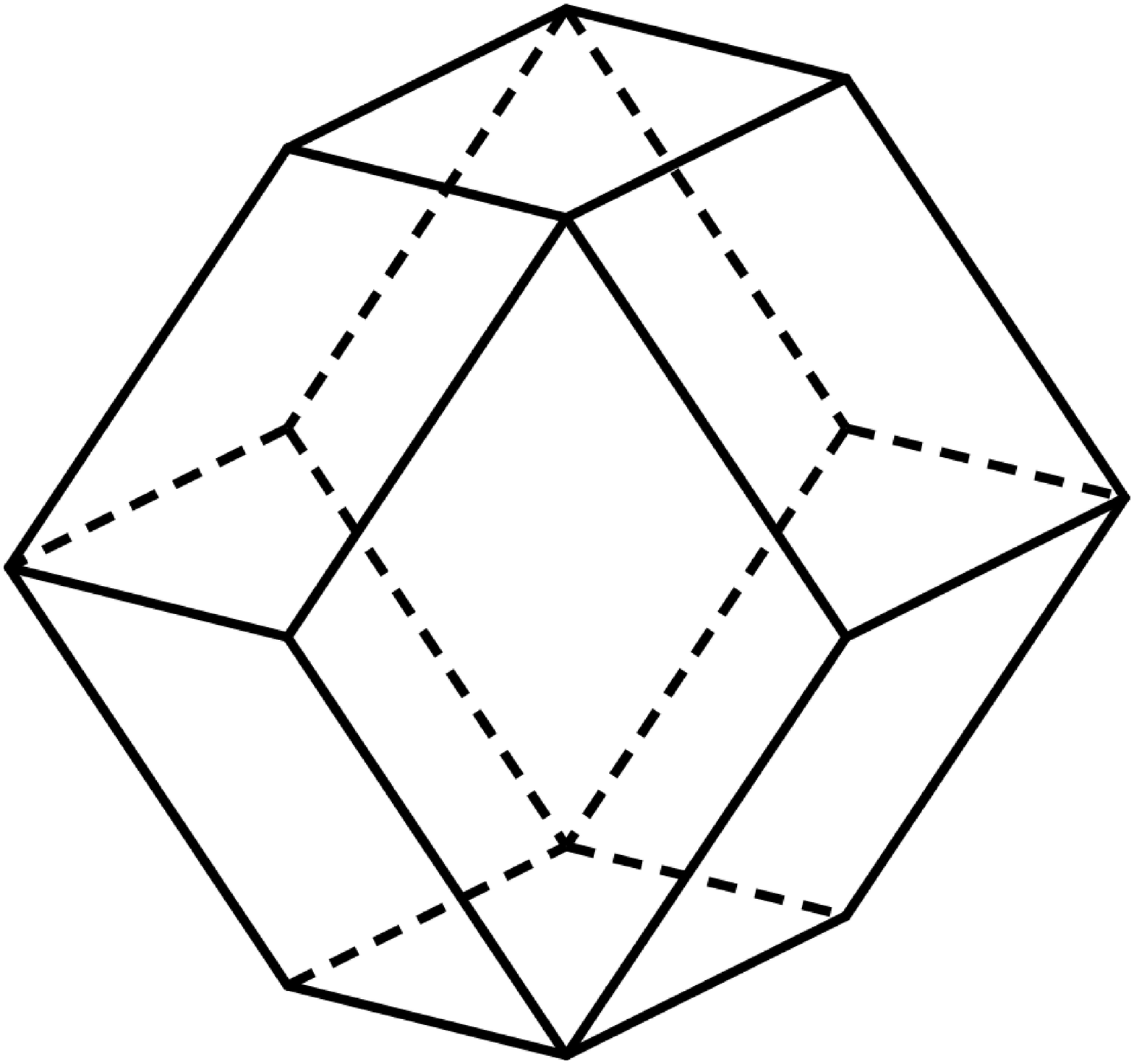}
\vspace*{-0,4cm}\caption{$\,$}
\label{f2}
\end{figure}}

\vspace*{.5cm}
{\bf Stratifications}
%(\cite{Dei}, Sect. 1 and Th. 5.11) :}

%We now mention without proofs some basic standard results from the Theory of Lie groups acting on finite dimensional manifolds. We refer the reader to \cite{Dei} and \cite{Kaw} for details.

%To determine the number $\dim \langle P\rangle$, we may fix $G=G(P)$ and confine attention to polyhedra $Q\cong P$ whose centroid is the same as that of $P$, say $0$ and whose symmetry group is $G$ itself. Thus $G$ is a subgroup of $\mathcal{O}(3)$.

For a transformation group $(M,G)$ the structure of the orbit space $M/G$ usually is complicated, for example it is not necessarily a manifold.

However, when the Lie group $G$ is compact and the manifold $M$ is without boundary it can be shown that they are stratified into smooth manifolds.

We now describe briefly a stratification of  $\mathbb{E}^3$ associated with a finite subgroup $G$ of the orthogonal group $O(3)$ that will help us to understand the relationship between the action of $G$ and the number $\dim \langle P\rangle $.

\begin{defi} ({\rm Stratification}) Let $X$ be a topological subspace of some Euclidian space  $\mathbb{E}^3$.

A partition  $\sum=\{M_i \, | \, i=1,\ldots,k\}$ of (pairwise disjoint) subsets of $X$ is called a stratification of $X$ if $\sum$ satisfies the followings:

\begin{enumerate}
\item Each $M_i,\,i=1,\ldots,k$ is a connected smooth submanifold of $\mathbb{E}^3$, called a \newline $\sum$-stratum.

\item For each $i$, the closure $\overline{M_i}$ is the union of $M_i$ and the $M_j$'s with lower dimensions than the dimension of $M_i$, that is, the relative closure $X \cap \overline{M_i}$ is the union of elements of $\sum$, one being $M_i$ itself and the others being of dimension less than the $\dim M_i$.

\end{enumerate}

This condition is called {\it frontier condition}. The dimension $\dim(X)$ is
\vspace*{-0.4cm}
$$\max \{\dim(M_i) \mid i=1,\ldots,k\} \ .$$
\end{defi}

The stratification mainly is done by the help of the theorem so called Slice Theorem which is fundamental in studding the structure of the transformation groups (see \cite{Kaw}, Th. 4.11).

\vspace*{.25cm}
Let $X=\mathbb{E}^3$ and $G$ a compact subgroup of $O(3)$, acting via $\phi$ on $\mathbb{E}^3$ as above.
%Our basic reference for this part are \cite{Kaw}, chapter 4, and \cite{Rob}, chapter 3.
For each  $x \in \mathbb{E}^3$ let $\: \: G_x = \{g \in G \mid g(x)=x\} \: \:$ be the isotropy subgroup of $G$ at $x$ and $F_x=Fix \left( G_x \right)$ be the set of all fixed points of $G_x$. Thus
\vspace*{-0,2cm}
$$F_x=\{y \in \mathbb{E}^3 \mid \text{ for all } g \in G_x, \, g(y)=y\} \, .$$

Since $x \in F_x$ and for all $g \in G$, $g(0)=0$  and for $y,z$ in $F_x$  and  $\lambda,\mu\in \mathbb{R}$ we have $g(\lambda y+\mu z)=\lambda y+\mu z$, $F_x$  is a linear subspace of $\mathbb{E}^3$. Define an equivalence relation $\sim_{G}$  on $\mathbb{E}^3$  as follows. Put $x\sim_{G} y$  if $F_x=F_y$.

Now let $x \in \mathbb{E}^3$ and $y \in F_x$. If $F_x=F_y$ then, $y \in [x]$, the equivalence class of $x$ in  $\sim_{G}$.  However $y \in F_x$  implies that $F_y \subseteq F_x$, since  $G_x \subseteq G_y$. Thus $ F_y$  is a linear subspace of  $F_x$. For $z\in \mathbb{E}^3-F_x$, we cannot have $z \sim_{G} x$. Therefore $[x]=\{y\in F_x: F_y=F_x\}$. So $[x]$ is complement in  $F_x$ of finitely many subspaces of $F_x$. Hence the equivalence classes $[x]$, $x \in \mathbb{E}^3$ stratify $\mathbb{E}^3$ with finitely many such strata (orbit types) (see \cite{Dei}, Th. 5.11). The dimension of $F_x$ denoted by $\delta(x)$ is called the \emph{degree of freedom} of $x$.

For example, let $G$ be a group generated by rotation matrix
%\vspace*{-.1cm}
$$A=\left(
    \begin{array}{ccc}
      \cos(\theta) & -\sin(\theta) & 0 \\
      \sin(\theta)& \cos(\theta) & 0 \\
      0 & 0 & 1 \\
    \end{array}
  \right),\,\,\,\,\,\,\,\theta=\frac{2\pi}{n},
$$

%\vspace*{-.3cm}
\noindent about $z$-axis through $\theta$.

Then, there are just four strata under this group action on $\mathbb{E}^3$, namely the $0$-stratum $\{0\}$, the open rays $x=y=0, z>0$ and $x=y=0, z<0$, and the complement of the $z$-axis.

Now let $P$ be a polyhedron with $G=G(P)$ some finite subgroup of $\mathcal{O}(3)$ and $Q\in \langle P \rangle$ as above.
Under the restriction imposed on $Q$ within the symmetry type of $P$, each vertex of $Q$ may be moved along a line, or within a plane, or in any direction in  $\mathbb{E}^3$ near (without changing the symmetry type) its initial position in $P$ itself, having one, two or three degrees of freedom. Likewise, each face $F$ of $Q$ may have one, two or three degrees of freedom close to the corresponding face of $P$, according as F intersects a 1-stratum in an interior point of $F$ (necessarily at right angles), or intersects a 2-stratum in the interior of $F$ (again at right angles), or neither of these.

Considering the action of $G=G(P)$ on $F(P)$, let
%\vspace*{-.3cm}
$$
\overline{v} \, = \, G(v) \, = \, \left\{ g(v) \, \mid \, g\in G \right\} \: , \: \: \overline{F} \, = \, G(F) \, = \, \left\{ g(F) \, \mid \, g\in G \right\} \, , \, \,\, \text{ and }
$$
\vspace*{-.6cm}
$$
F_0(P)/ G \, = \, \left\{ \overline{v} \, \mid \, v\in F_0(P) \right\} \: , \: \: F_2(P)/ G \, = \, \left\{ \overline{F} \, \mid \, F\in F_2(P) \right\}
$$

%\vspace*{-.2cm}
\noindent be the collection of orbits under the action of $G$ on face lattice $F(P)$.

%Note that (Definition \ref{TwoP})
Define $\mathcal{M}_P$ to be the set of all pairs $(v,F)\in F_0(P) \times F_2(P)$ for which $v\in F$ and
\vspace*{-.2cm}
$$
\mathcal{M}_P/ G \, = \, \left\{ \left( g(v), g(F) \right) \, \mid \, g\in G \text{ and } (v,F)\in \mathcal{M}_P \right\} \, .
$$

%\vspace*{-.3cm}
Clearly $\mu (P)$ and $\mu _{*}(P)$ are the cardinalities of $\mathcal{M}_P$ and $\mathcal{M}_P/ G$, respectively.

\vspace*{.2cm}
We observe that if $v,u\in F_0(P)$ and $\overline{v}=\overline{u}$ then,
\vspace*{-.2cm}
$$
\dim(v) \, = \, \dim(u) \, .
$$

\vspace*{-.2cm}
Therefore we can define the fixed dimensions $\delta (\xi)$, $\xi \in F_0(P)/ G$, $\delta (\xi)=\dim(v)$, the degree of freedom of $v$ for an arbitrary $\overline{v}\in \xi$.
The same holds for $\delta (\zeta)$, $\zeta \in F_2(P)/ G$.

Our aim is to count the number of the vertices and the faces with dimensions $k=1,2,3$ and then, by subtracting the independent incidences $\mu _{*}(P)$, express the $\dim \langle P \rangle$ in terms of edge orbits alone.

%, that is
%\vspace*{-.3cm}
%$$\dim\langle P \rangle= \left( \sum (\delta(\overline{v}) + \sum \delta(\overline{F}))\right)-\mu_*(P) \, ,$$
%
%\vspace*{-.3cm}
%\noindent $\sum$ is extended over all vertex orbits and face orbits of $P$, respectively.

%Let $\delta(v)=k$, $k=1,2,3$, if $v$ lies on a $k$-stratum of $G$, and let  $\delta(F)=k$ with $k=1,2,3$, if the face $F$ has $k$ degrees of freedom with the above restrictions. Our aim is to consider the action $G(P)$ on vertices, edges and faces of $P$, by counting the number of the vertices $v$ with $\delta(v)=k$ and the number of faces  $\delta(F)$ with  $\delta(F)=k, k=1,2,3$ to express the $\dim\langle P \rangle$ in terms of edge orbits alone.

%Define $[v]=\{ g(v):g \in G\}$  and  $[F]=\{ g(F):g \in G\}$. Let $F_0(P)/G=\{[v]: v \in F_0(P)\}$  and  $F_2(P)/G=\{[F]: F \in F_2(P)\}$ be the collection of all vertex orbits and face orbits of $P$ respectively.

\section{Fundamental regions and main theorem}
%\section{Proof of Main Theorem}

In this Section we consider finite subgroups of isometries which are generated by reflections namely $[q], [2, q], [3, 3], [3, 4]$ and $[3,5]$ (Table 2).

\vspace*{.2cm}
The finite subgroups of $\mathbb{E}(3)$ which are generated by reflections in the plane, are given in the following table. They are called reflection groups, for obvious reason.
\begin{center}
\vspace*{-.2cm} {\small{\begin{tabular}{|c|l|c|}
  \hline
  Symbol & description & order \\
  \hline
                  & $q=1$: One plane of reflection; $q\geq 2$: $q$ equally inclined planes  &   \\
    $[q],q\geq 1$ & of reflection passing through a $q$-fold axis of rotation, dihedral \hspace{11cm}& $2q$ \\
                  & reflection group. &   \\
 \hline
                  & $q$ equally inclined planes of reflection passing through a $q$-fold &   \\
    $[2,q]$       & axis of rotation and reflection in a equatorial plane.   & 4q \\
                  & $q$ 2-fold axes of rotation. The group of $q$-prism.&   \\
\hline
     $[3,3]$      & Four 3-fold and three 2-fold axes. Six planes of reflection. &  \\
                  & Symmetry group of the regular tetrahedron. & 24  \\
\hline
                  & Three 4-fold and four 3-fold and six 2-fold axes of rotation. &   \\
     $[3,4]$      & Nine planes of reflection. Symmetry group of the cube. & 48 \\
\hline
                  & Six 5-fold, ten 3-fold and fifteen 2-fold axes of rotation.  &   \\
     $[3,5]$      & Fifteen planes of reflection. Symmetry group of the icosahedron. & 120 \\
  \hline
\end{tabular}}}
\end{center}

%\vspace*{0,15cm}
\begin{center} Table 2: Reflection groups \end{center}

It is well known that the \emph{fundamental region} $\Delta$ for the action of $[3, 3], [3, 4], [3, 5]$ and $[2, q]$ on the sphere $S^2$  are spherical triangles [1]. For $[q]$ the dihedral reflection group generated by two reflections, the fundamental region is a ``lune'' of angle $\frac{\pi}{q}$.
We may use the fundamental region of a reflection group to construct a stratification of $\mathbb{E}^3$. For instance consider the tetrahedron $OABC$ (or its spherical projection) as a fundamental region of [3, 4] in Figure \ref{f3} (\cite{Rob}, p. 81).
We take the origin $O$ as a $0{\text{-stratum}}$. By removing the origin from the rays $OA$, $OB$ and $OC$ we get three $1{\text{-strata}}$, the interiors of the region $AOB$, $AOC$ and $BOC$ are $0{\text{-strata}}$.

Finally, the interior points of $\mathbb{E}^3$  bounded by sectors $AOB$, $AOC$ and $BOC$ is $3{\text{-strata}}$. By transferring these strata under the action of [3, 4] we obtain the required stratification of $\mathbb{E}^3$.

For other reflection groups a stratification of $\mathbb{E}^3$ is constructed in analogues fashion.
Returning to our main problem, we now consider the following notion. Suppose that a finite reflection group $G$ in $\mathcal{O}(3)$ has its fundamental region a spherical triangle $\Delta$ and, let $P$ be a polyhedron with $G(P)=G$.
%\vspace*{-.5cm}
%UMA FIGURA
{ \begin{figure}[!htb] \centering
\hspace*{.5cm} \includegraphics[scale=.2]{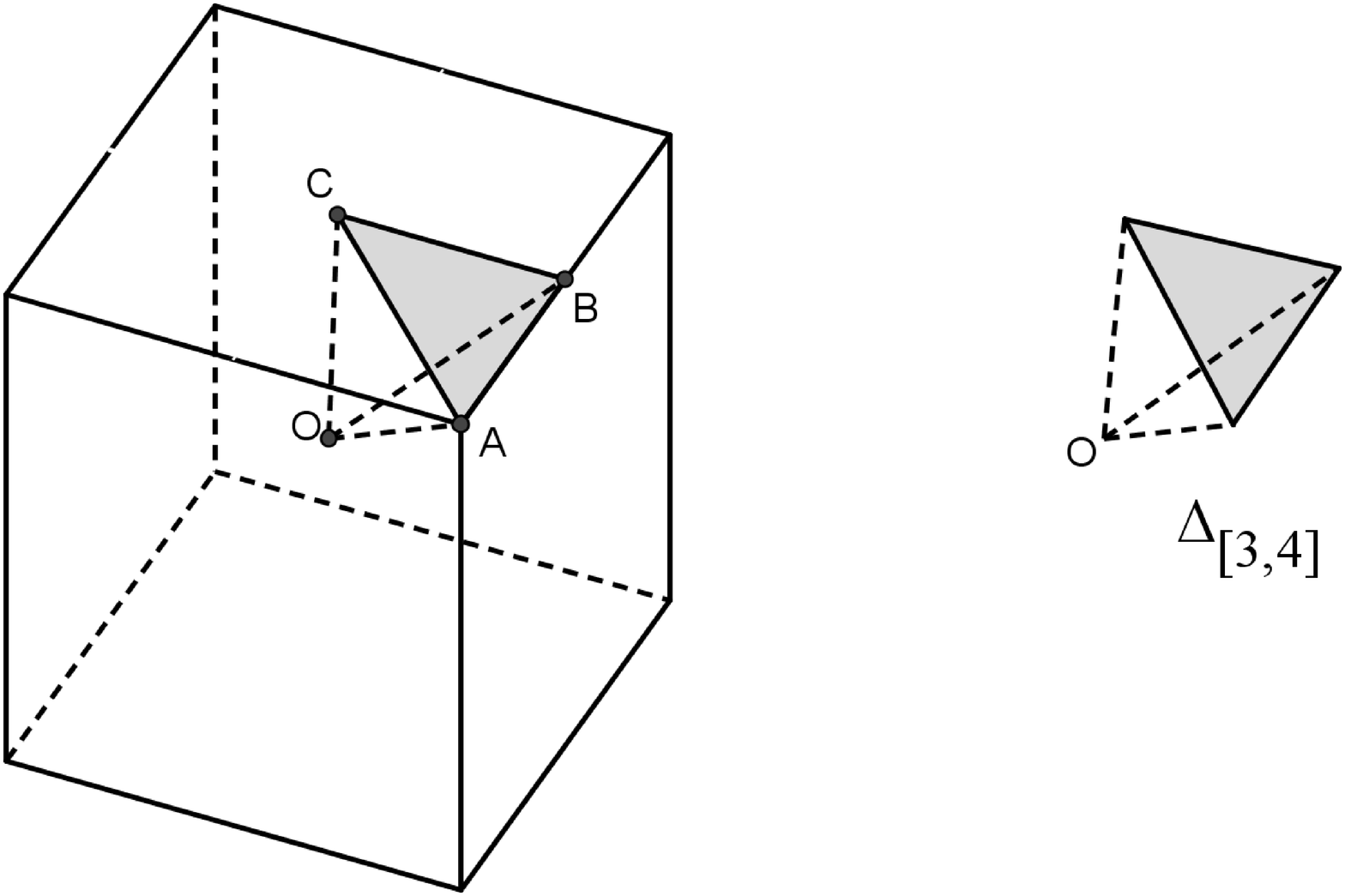}
\vspace*{-0,4cm}\caption{$\,$}
\label{f3}
\end{figure}}

We denote by  $\Delta_p$ that portion of the surface of $P$ (namely those vertices, edges and subpolygonal faces) which lie within $\Delta$, and call $\Delta_p$ a \emph{basic region} of $P$ (Figure \ref{f4}).
%\vspace*{-.2cm}
%UMA FIGURA
{ \begin{figure}[!htb] \centering
\hspace*{.75cm} \includegraphics[scale=.225]{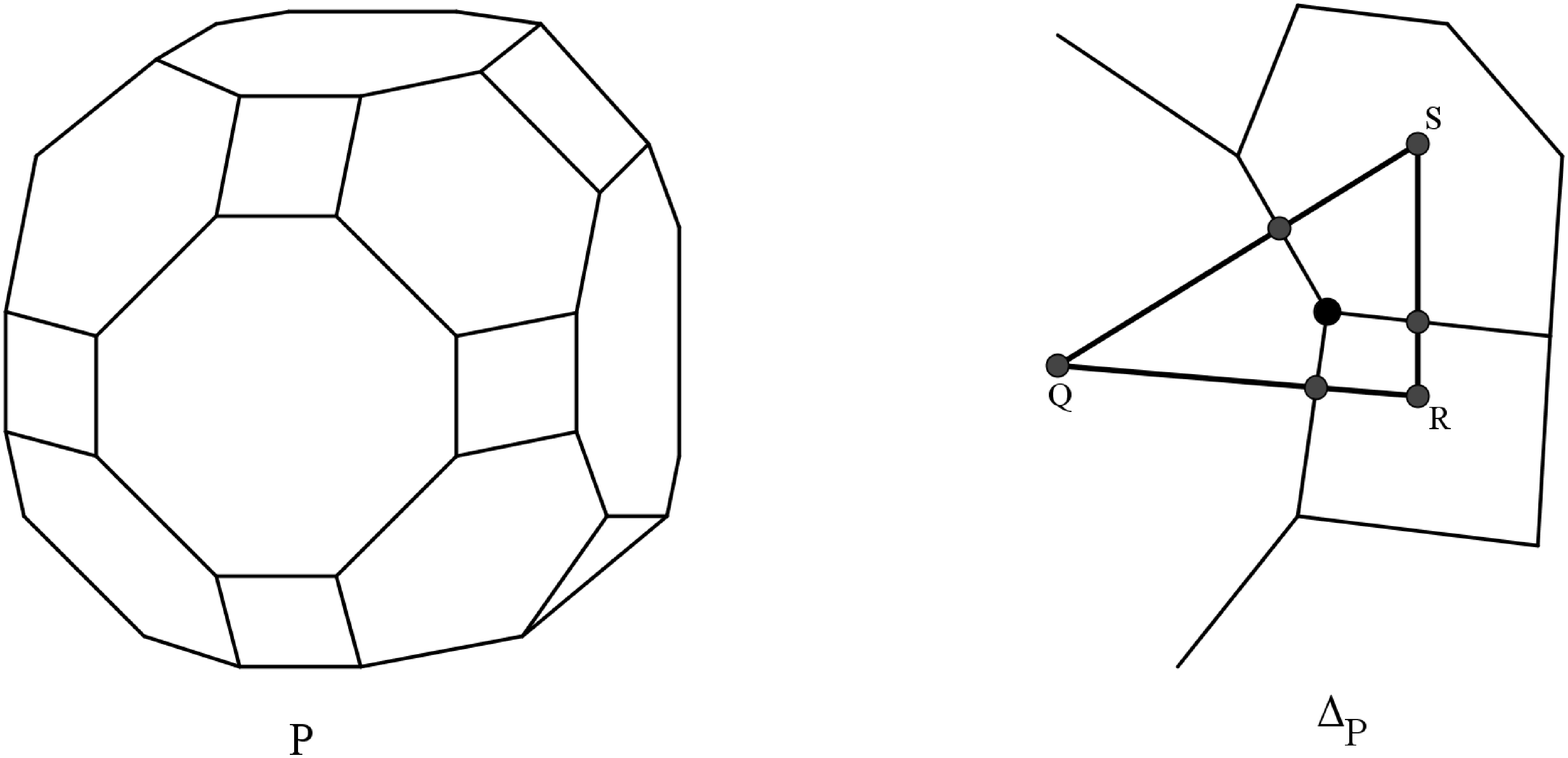}
\vspace*{-0,4cm}\caption{$\,$}
\label{f4}
\end{figure}}

\vspace*{-0.1cm}
Hence $\Delta_p$ is a simple closed planar polygonal region. If $n_{\Delta_p}, e_{\Delta_p}$ and $f_{\Delta_p}$ denote the total number of distinct vertices, edges and subpolygonal faces of  $\Delta_p$ respectively then, from Euler's formula by stereographic projection we get
\vspace*{-0,2cm}
$$n_{\Delta_p} - e_{\Delta_p} + f_{\Delta_p}=1 \, .$$

As illustrated example, let $P$ be truncated cuboctahedron (Figure \ref{f4}) with symmetry group $G(P)=[3,4]$ of cube, with order 48. Thus the basic region $\Delta_p$ has $n_{\Delta_p}=7$, $e_{\Delta_p}=9$, $f_{\Delta_p}=3$, and $n_{\Delta_p} - e_{\Delta_p} + f_{\Delta_p}=7-9+3=1$.

\vspace*{0.1cm}
Now having our necessary tools, we are in the position to state and prove our main theorem in this Section.

\begin{teo}Let $G$ be a finite reflection group in $\mathbb{E}(3)$ and $P$ a polyhedron with $G(P)=G$. Then, $\dim \langle P\rangle= \epsilon-1$ , where $\epsilon$ is the number of edge orbits of the action of $G$ on the set of edges of $P$.
\end{teo}

First we prove the following lemma.

\begin{lema} Assuming the hypothesis of the theorem, let  $\Delta_p$  be a basic region for $P$ such that the corners of fundamental region of $\Delta$ of $G$ are vertices of $P$ (see Figure \ref{f5}). Then, the number of incident pairs of vertices and faces of  $\Delta_p$, the multiplicity $\mu(\Delta_p)$  of  $\Delta_p$  is given by  $\mu(\Delta_p)=2e-\beta$ where $e$ is the total number of edges of   $\Delta_p$ and $\beta$ the number of vertices on the boundary of the fundamental region.
\end{lema}
{\bf Proof:} By adjoining an extra face, say $K$, to  $\Delta_p$, namely the complement of   $\Delta_p$ itself with respect to the sphere we get a map $M_P$ on sphere. But the number of edges (and vertices) of $M_P$ is equal to the number of edges (and vertices) of  $\Delta_p$. Hence  $\mu(\Delta_p)=2e$. Since there are $\beta$ vertices on boundary of  $\Delta$, with respect to that extra face $K$, we have  $\mu(\Delta_p)=\mu(P)-\beta=2e - \beta$.
\vspace*{-0.3cm}
%UMA FIGURA
{ \begin{figure}[!htb] \centering
\hspace*{.75cm} \includegraphics[scale=.15]{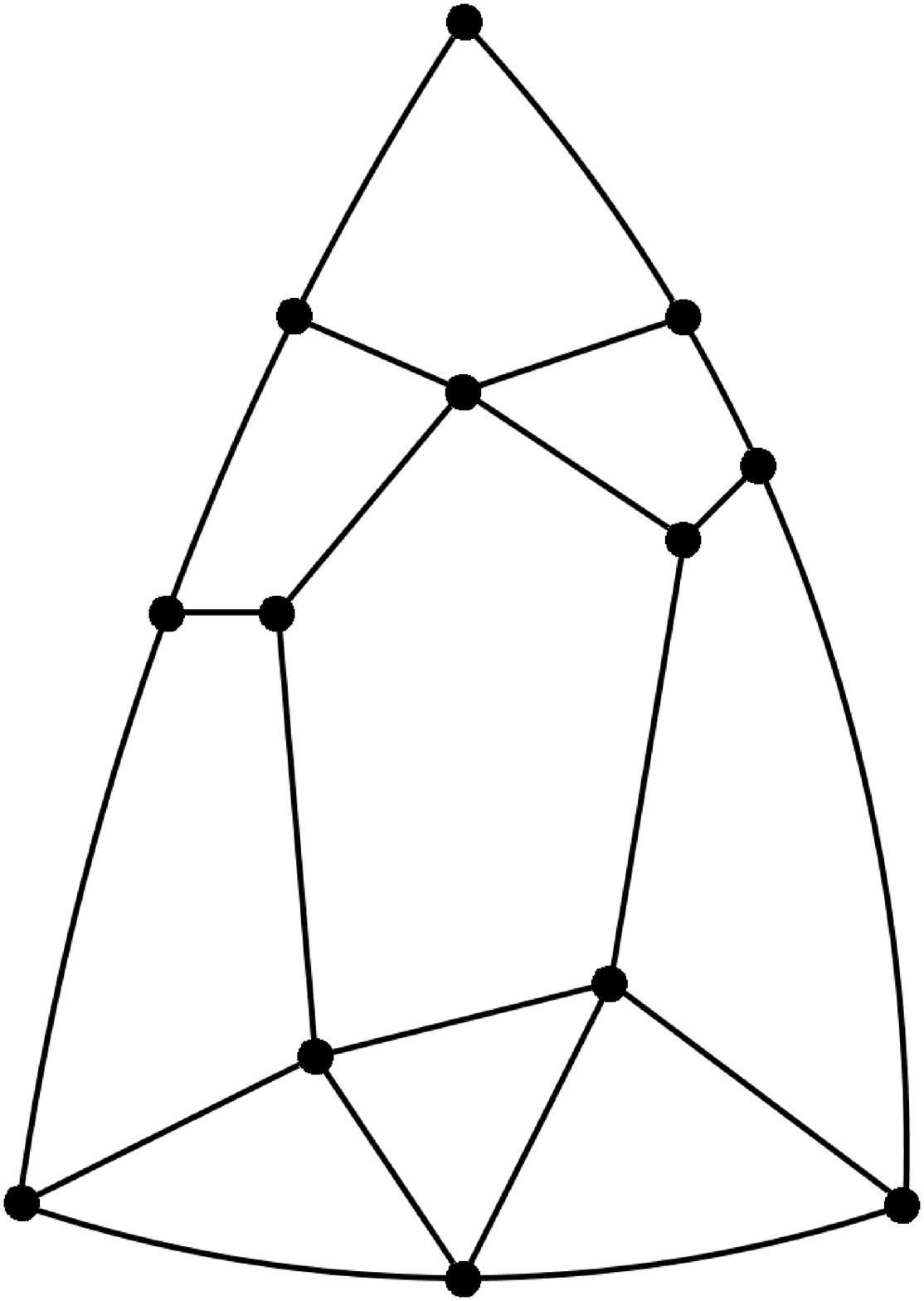}
\vspace*{-0,1cm}\caption{$\,$}
\label{f5}
\end{figure}}

\vspace*{-0,5cm}
$\:$ \hfill$\square$

\vspace*{-0,1cm}
 \begin{obs}
  {\rm Let $P$ be a polyhedron and $v$ a vertex of $P$. According to stratification of  $\mathbb{E}^3$ with reflection group $G(P)$ of $P$, $\delta(v)$, the degree of freedom of $v$,  is one, two or three if $v$ is on a corner or side or within the interior of   $\Delta_p$ respectively. Similarly if a face $F$ of $P$ has as its interior point a corner of $\Delta$ then, $\delta(F)=1$. If $F$ is orthogonal to a side of  $\Delta_p$  or lies inside  $\Delta_p$  then, $\delta(F)=2 $ or 3, respectively.}
 \end{obs}

 \noindent {\bf Proof of theorem:} Let  $\eta(1),\eta(2)$ and $\eta(3)$ be the number of vertices of $\Delta_p$   and   $\phi(1),\phi(2)$ and $\phi(3)$ the number of faces of $\Delta_p$ with one, two and three degrees of freedom respectively.

 First we assume that $P$ has no face $F$ with $\delta(F)$ equal to one or two. Then, $\mu_*(P)=\mu(\Delta_p)=2e-\beta$. After factoring out the effect of dilation we get
 $$
 \begin{array}{rcl}
   \dim\langle P \rangle & = & 1\eta(1)+2\eta(2)+3\eta(3)+3\phi(3)-\mu(\Delta_p)-1 \\[8pt]
     & = & 3(n_{\Delta_p}+\phi_{\Delta_p})-2e_{\Delta_p}-4 \, ,
 \end{array}$$
\noindent since $\eta(1)=3$.

%-------- ILDA
But $n_{\Delta_p}+\phi_{\Delta_p}=e_{\Delta_p}+1$. Hence,
$$\dim\langle P \rangle=3(e_{\Delta_p}+1)-2e_{\Delta_p}-4= e_{\Delta_p}-1 \, .$$
\noindent Now clearly the number of edges of  $\Delta_p$  is exactly the number of edge orbits of $P$. Therefore the theorem follows in this case.

Next suppose $P$ has one face $F$ with $\delta(F)=1$. This means that there is a face $F$ such that one corner say $v$ of $\Delta_p$   is an interior point of $F$.

Let $P'$ be a polyhedron which we get, by changing ``fake'' edges in $\Delta_p$ into real ones (see Figure \ref{f6}). This is done as follows.

We remove the constraint that the plane of  $F \cap \Delta_p$  is perpendicular to the ray $ov$. The edges of  $F \cap \Delta_p$  that lie in the boundary of  $\Delta_p$  are also edges of $P'$ where $P'$ has a basic region  $\Delta_{p'}$ , say, with the same combinatorial structure as $\Delta_p$, and with vertices arbitrarily close to those of $\Delta_p$.
Since we have substituted a face with one degree of freedom by a face with three degrees of freedom and since the new vertex and its incidence cancel each other and hence do not effect our calculation for  $ \dim\langle P \rangle $, we have $ \dim\langle P' \rangle = \dim\langle P \rangle  + 2$. But $\dim\langle P' \rangle = \epsilon'-1$ and $\epsilon'=\epsilon+2$ with obvious notations. Hence
$$\dim\langle P \rangle = \dim\langle P' \rangle -2= \epsilon'-1-2=\epsilon-1 \, .$$

The process will continue if $P$ has two or three (on other corners of $\Delta$) faces of degree one.
%UMA FIGURA
{ \begin{figure}[!htb] \centering
\vspace*{-.5cm} \hspace*{.25cm} \includegraphics[scale=.175]{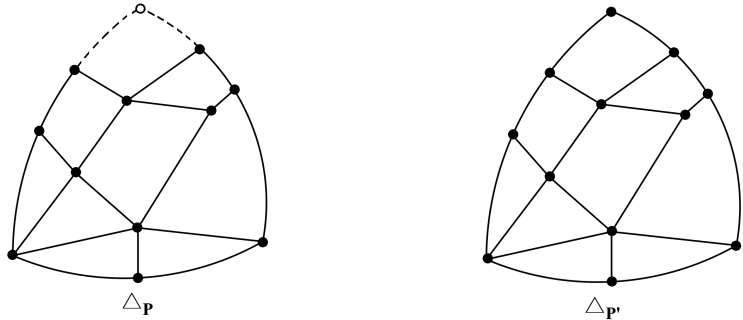}
\vspace*{-0,4cm}\caption{\small{\emph{Broken lines represent ``fake'' edges and ``o'' a fake vertex.}}}
\label{f6}
\end{figure}}
%\vspace*{-.5cm}
%\begin{center}
%{\small{Figure 6}\footnote{Broken lines represent "fake" edges and "o"a fake vertex.}}
%\end{center}

Here we remark that, in the process of changing ``fake'' edges in $\Delta_p$  into real ones, since the transforms of $\Delta_{p'}$ under corresponding group action is a polyhedral graph (planar and 3-connected), by the Theorem of Steinitz \cite{Stei}, there exists a polyhedron P' which geometrically realizes $\Delta_{P'}$.

\vspace*{.25cm}
Finally, suppose $P$ has a face $F$ with $\delta(F)=2$ (Figure \ref{f7}). We construct $P'$ by adjoining the fake edge to $\Delta_p$. Then, $\dim\langle P' \rangle $ differs by one from
$\dim\langle P \rangle$ for replacement of a face with two degrees of freedom, by a face of three degrees of freedom. Hence
$\dim\langle P' \rangle = \dim\langle P \rangle +1$. Because of $\dim\langle P' \rangle =\epsilon'-1$ and $\epsilon=\epsilon'-1$, we have
\vspace*{-0,3cm}
$$\dim\langle P \rangle +1=\epsilon'-1=\epsilon\,\,\,\,\text{ and }\,\,\,\, \dim\langle P \rangle =\epsilon-1 \, .$$

%UMA FIGURA
{ \begin{figure}[!htb] \centering
\hspace*{.25cm} \includegraphics[scale=.175]{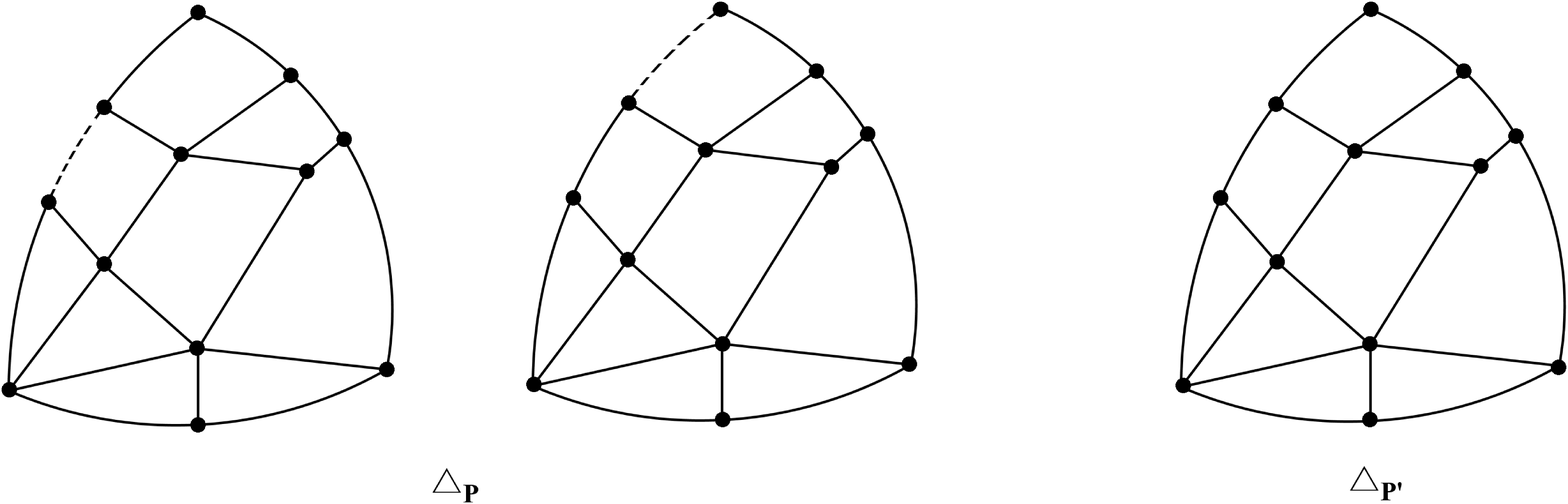}
\vspace*{-0,4cm}\caption{$\,$}
\label{f7}
\end{figure}}

\newpage
Now this inductive process can be continued if $P$ has any number of faces with two degrees of freedom. In each step of construction, $\dim\langle P \rangle$ and $\epsilon$ each increase by one, while the operation leaves every other quantity in our calculation fixed. For the case of the group $[q]$ where $\Delta$ is a ``lune'', the proof proceeds in similar way and is omitted here to avoid repetition. The proof of the theorem now is complete. \hfill$\square$

\end{document}